\documentclass[12pt]{article}
\usepackage{amssymb}
\usepackage{amsfonts}
\usepackage{amsmath}
\usepackage{amsbsy}
\newcommand{\be}{\begin{equation}}
\newcommand{\ee}{\end{equation}}
\newcommand{\bea}{\begin{eqnarray}}
\newcommand{\eea}{\end{eqnarray}}
\newcommand{\ba}{\begin{array}}
\newcommand{\ea}{\end{array}}

\newcommand{\bc}{\begin{center}}
\newcommand{\ec}{\end{center}}
\newcommand{\ben}{\begin{enumerate}}
\newcommand{\een}{\end{enumerate}}
\newcommand{\bfi}{\begin{figure}}
\newcommand{\efi}{\end{figure}}

\newcommand{\bq}{\begin{quote}}
\newcommand{\eq}{\end{quote}}
\newcommand{\bqu}{\begin{quotation}}
\newcommand{\equ}{\end{quotation}}
\newenvironment{emphit}{\begin{itemize}}{\end{itemize}}
\newcommand{\bemp}{\begin{emphit}}
\newcommand{\eemp}{\end{emphit}}

\newcommand{\bt}{\begin{tabular}}
\newcommand{\et}{\end{tabular}}

\newtheorem{myth}{Theorem}[section]
\newtheorem{mylem}{Lemma}[section]
\newtheorem{mycor}{Corollary}[section]

\newtheorem{mydef}{Definition}[section]
\newtheorem{myrem}{Remark}[section]
\newtheorem{myquest}{Question}[section]

\begin{document}
\date{}
\title{Basic Properties
Of Second Smarandache Bol Loops 
}
\author{T\`em\'it\d{\'o}p\d{\'e} Gb\d{\'o}l\'ah\`an Ja\'iy\'e\d ol\'a\thanks{All correspondence to be addressed to this author}\\
Department of Mathematics,\\
Obafemi Awolowo University, Ile Ife, Nigeria.\\
jaiyeolatemitope@yahoo.com, tjayeola@oauife.edu.ng} \maketitle
\begin{abstract}
The pair $(G_H,\cdot )$ is called a special loop if $(G,\cdot )$ is
a loop with an arbitrary subloop $(H,\cdot )$. A special loop
$(G_H,\cdot )$ is called a second Smarandache Bol
loop(S$_{2^{{\tiny\textrm{nd}}}}$BL) if and only if it obeys the
second Smarandache Bol identity $(xs\cdot z)s=x(sz\cdot s)$ for all
$x,z$ in $G$ and $s$ in $H$. The popularly known and well studied
class of loops called Bol loops fall into this class and so
S$_{2^{{\tiny\textrm{nd}}}}$BLs generalize Bol loops. The basic
properties of S$_{2^{{\tiny\textrm{nd}}}}$BLs are studied. These
properties are all Smarandache in nature. The results in this work
generalize the basic properties of Bol loops, found in the Ph.D.
thesis of D. A. Robinson. Some questions for further studies are
raised.
\end{abstract}

\section{Introduction}
The study of the Smarandache concept in groupoids was initiated by
W. B. Vasantha Kandasamy in \cite{phd86}. In her book \cite{phd75}
and first paper \cite{phd83} on Smarandache concept in loops, she
defined a Smarandache loop(S-loop) as a loop with at least a subloop
which forms a subgroup under the binary operation of the loop. The
present author has contributed to the study of S-quasigroups and
S-loops in \cite{sma1,sma2,sma3,sma14,sma15,sma8,sma9,sma16} by introducing some new
concepts immediately after the works of Muktibodh \cite{muk1,muk2}. His recent monograph \cite{smabook1} gives inter-relationships and connections between and
among the various Smarandache concepts and notions that have been
developed in the aforementioned papers.
  
But in the quest of developing the concept of Smarandache
quasigroups and loops into a theory of its own just as in
quasigroups and loop theory(see
\cite{phd41,phd39,phd49,phd42,phd3,phd75}), there is the need to
introduce identities for types and varieties of Smarandache
quasigroups and loops. For now, a Smarandache loop or Smarandache
quasigroup will be called a first Smarandache
loop(S$_{1^{{\tiny\textrm{st}}}}$-loop) or first Smarandache
quasigroup(S$_{1^{{\tiny\textrm{st}}}}$-quasigroup).

Let $L$ be a non-empty set. Define a binary operation ($\cdot $) on
$L$ : if $x\cdot y\in L$ for all $x,y\in L$, $(L, \cdot )$ is called
a groupoid. If the system of equations ; $a\cdot x=b$ and $y\cdot
a=b$ have unique solutions for $x$ and $y$ respectively, then $(L,
\cdot )$ is called a quasigroup. For each $x\in L$, the elements
$x^\rho =xJ_\rho ,x^\lambda =xJ_\lambda\in L$ such that
$xx^\rho=e^\rho$ and $x^\lambda x=e^\lambda$ are called the right,
left inverses of $x$ respectively. Furthermore, if there exists a
unique element $e=e_\rho =e_\lambda$ in $L$ called the identity
element such that for all $x$ in $L$, $x\cdot e=e\cdot x=x$, $(L,
\cdot )$ is called a loop. We write $xy$ instead of $x\cdot y$, and
stipulate that $\cdot$ has lower priority than juxtaposition among
factors to be multiplied. For instance, $x\cdot yz$ stands for
$x(yz)$. A loop is called a right Bol loop(Bol loop in short) if and
only if it obeys the identity
\begin{displaymath}
(xy\cdot z)y=x(yz\cdot y).
\end{displaymath}
This class of loops was the first to catch the attention of loop
theorists and the first comprehensive study of this class of loops
was carried out by Robinson \cite{phd85}.

The aim of this work is to introduce and study the basic properties
of a new class of loops called second Smarandache Bol
loops(S$_{2^{{\tiny\textrm{nd}}}}$BLs). The popularly known and well
studied class of loops called Bol loops fall into this class and so
S$_{2^{{\tiny\textrm{nd}}}}$BLs generalize Bol loops. The basic
properties of S$_{2^{{\tiny\textrm{nd}}}}$BLs are studied. These
properties are all Smarandache in nature. The results in this work
generalize the basic properties of Bol loops, found in the Ph.D.
thesis \cite{phd85} and the paper \cite{phd194} of D. A. Robinson.
Some questions for further studies are raised.

\section{Preliminaries}
\begin{mydef}\label{1:1}
Let $(G,\cdot )$ be a quasigroup with an arbitrary non-trivial
subquasigroup $(H,\cdot )$. Then, $(G_H,\cdot )$ is called a special
quasigroup with special subquasigroup $(H,\cdot )$. If $(G,\cdot )$
is a loop with an arbitrary non-trivial subloop $(H,\cdot )$. Then,
$(G_H,\cdot )$ is called a special loop with special subloop
$(H,\cdot )$. If $(H,\cdot )$ is of exponent $2$, then $(G_H,\cdot
)$ is called a special loop of Smarandache exponent $2$.

A special quasigroup $(G_H,\cdot )$ is called a second Smarandache
right Bol quasigroup(S$_{2^{{\tiny\textrm{nd}}}}$-right Bol
quasigroup) or simply a second Smarandache Bol
quasigroup(S$_{2^{{\tiny\textrm{nd}}}}$-Bol quasigroup) and
abbreviated S$_{2^{{\tiny\textrm{nd}}}}$RBQ or
S$_{2^{{\tiny\textrm{nd}}}}$BQ if and only if it obeys the second
Smarandache Bol identity(S$_{2^{{\tiny\textrm{nd}}}}$-Bol identity)
i.e S$_{2^{{\tiny\textrm{nd}}}}$BI
\begin{equation}\label{eq:1}
(xs\cdot z)s=x(sz\cdot s)~\textrm{for all}~x,z\in
G~\textrm{and}~s\in H.
\end{equation}
Hence, if $(G_H,\cdot )$ is a special loop, and it obeys the
S$_{2^{{\tiny\textrm{nd}}}}$BI, it is called a second Smarandache
Bol loop(S$_{2^{{\tiny\textrm{nd}}}}$-Bol loop) and abbreviated
S$_{2^{{\tiny\textrm{nd}}}}$BL.
\end{mydef}
\begin{myrem}
A Smarandache Bol loop(i.e a loop with at least a non-trivial
subloop that is a Bol loop) will now be called a first Smarandache
Bol loop(S$_{1^{{\tiny\textrm{st}}}}$-Bol loop). It is easy to see
that a S$_{2^{{\tiny\textrm{nd}}}}$BL is a
S$_{1^{{\tiny\textrm{nd}}}}$BL. But the reverse is not generally
true. So S$_{2^{{\tiny\textrm{nd}}}}$BLs are particular types of
S$_{1^{{\tiny\textrm{nd}}}}$BL. There study can be used to
generalise existing results in the theory of Bol loops by simply
forcing $H$ to be equal to $G$.
\end{myrem}

\begin{mydef}\label{1:2}
Let $(G,\cdot )$ be a quasigroup(loop). It is called a right inverse
property quasigroup(loop)[RIPQ(RIPL)] if and only if it obeys the
right inverse property(RIP) $yx\cdot x^\rho=y$ for all $x,y\in G$.
Similarly, it is called a left inverse property
quasigroup(loop)[LIPQ(LIPL)] if and only if it obeys the left
inverse property(LIP) $x^\lambda\cdot xy=y$ for all $x,y\in G$.
Hence, it is called an inverse property quasigroup(loop)[IPQ(IPL)]
if and only if it obeys both the RIP and LIP.

$(G,\cdot )$ is called a right alternative property
quasigroup(loop)[RAPQ(RAPL)] if and only if it obeys the right
alternative property(RAP) $y\cdot xx=yx\cdot x$ for all $x,y\in G$.
Similarly, it is called a left alternative property
quasigroup(loop)[LAPQ(LAPL)] if and only if it obeys the left
alternative property(LAP) $xx\cdot y=x\cdot xy$ for all $x,y\in G$.
Hence, it is called an alternative property
quasigroup(loop)[APQ(APL)] if and only if it obeys both the RAP and
LAP.

The bijection $L_x : G\rightarrow G$ defined as $yL_x=x\cdot y$ for
all $x,y\in G$ is called a left translation(multiplication) of $G$
while the bijection $R_x : G\rightarrow G$ defined as $yR_x=y\cdot
x$ for all $x,y\in G$ is called a right translation(multiplication)
of $G$.

$(G,\cdot )$ is said to be a right power alternative property
loop(RPAPL) if and only if it obeys the right power alternative
property(RPAP)
\begin{displaymath}
xy^n=\underbrace{(((xy)y)y)y\cdots
y}_{\textrm{$n$-times}}~\textrm{i.e.}~R_{y^n}=R_y^n~\textrm{for
all}~x,y\in G~\textrm{and}~n\in \mathbb{Z}.
\end{displaymath}

The right nucleus of $G$ denoted by $N_\rho (G, \cdot )=N_\rho
(G)=\{a\in G : y\cdot xa=yx\cdot a~\forall~ x, y\in G\}$.

Let $(G_H,\cdot )$ be a special quasigroup(loop). It is called a
second Smarandache right inverse property
quasigroup(loop)[S$_{2^{{\tiny\textrm{nd}}}}$RIPQ(S$_{2^{{\tiny\textrm{nd}}}}$RIPL)]
if and only if it obeys the second Smarandache right inverse
property(S$_{2^{{\tiny\textrm{nd}}}}$RIP) $ys\cdot s^\rho=y$ for all
$y\in G$ and $s\in H$. Similarly, it is called a second Smarandache
left inverse property
quasigroup(loop)[S$_{2^{{\tiny\textrm{nd}}}}$LIPQ(S$_{2^{{\tiny\textrm{nd}}}}$LIPL)]
if and only if it obeys the second Smarandache left inverse
property(S$_{2^{{\tiny\textrm{nd}}}}$LIP) $s^\lambda\cdot sy=y$ for
all $y\in G$ and $s\in H$. Hence, it is called a second Smarandache
inverse property
quasigroup(loop)[S$_{2^{{\tiny\textrm{nd}}}}$IPQ(S$_{2^{{\tiny\textrm{nd}}}}$IPL)]
if and only if it obeys both the S$_{2^{{\tiny\textrm{nd}}}}$RIP and
S$_{2^{{\tiny\textrm{nd}}}}$LIP.

$(G_H,\cdot )$ is called a third Smarandache right inverse property
quasigroup(loop)[S$_{3^{{\tiny\textrm{rd}}}}$RIPQ(S$_{3^{{\tiny\textrm{rd}}}}$RIPL)]
if and only if it obeys the third Smarandache right inverse
property(S$_{3^{{\tiny\textrm{rd}}}}$RIP) $sy\cdot y^\rho=s$ for all
$y\in G$ and $s\in H$.

$(G_H,\cdot )$ is called a second Smarandache right alternative
property
quasigroup(loop)[S$_{2^{{\tiny\textrm{nd}}}}$RAPQ(S$_{2^{{\tiny\textrm{nd}}}}$RAPL)]
if and only if it obeys the second Smarandache right alternative
property(S$_{2^{{\tiny\textrm{nd}}}}$RAP) $y\cdot ss=ys\cdot s$ for
all $y\in G$ and $s\in H$. Similarly, it is called a second
Smarandache left alternative property
quasigroup(loop)[S$_{2^{{\tiny\textrm{nd}}}}$LAPQ(S$_{2^{{\tiny\textrm{nd}}}}$LAPL)]
if and only if it obeys the second Smarandache left alternative
property(S$_{2^{{\tiny\textrm{nd}}}}$LAP) $ss\cdot y=s\cdot sy$ for
all $y\in G$ and $s\in H$. Hence, it is called an second Smarandache
alternative property
quasigroup(loop)[S$_{2^{{\tiny\textrm{nd}}}}$APQ(S$_{2^{{\tiny\textrm{nd}}}}$APL)]
if and only if it obeys both the S$_{2^{{\tiny\textrm{nd}}}}$RAP and
S$_{2^{{\tiny\textrm{nd}}}}$LAP.

$(G_H,\cdot )$ is said to be a Smarandache right power alternative
property loop(SRPAPL) if and only if it obeys the Smarandache right
power alternative property(SRPAP)
\begin{displaymath}
xs^n=\underbrace{(((xs)s)s)s\cdots
s}_{\textrm{$n$-times}}~\textrm{i.e.}~R_{s^n}=R_s^n~\textrm{for
all}~x\in G,~s\in H~\textrm{and}~n\in \mathbb{Z}.
\end{displaymath}

The Smarandache right nucleus of $G_H$ denoted by $SN_\rho (G_H,
\cdot )=SN_\rho (G_H)=N_\rho (G)\cap H$. $G_H$ is called a
Smarandache right nuclear square special loop if and only if $s^2\in
SN_\rho (G_H)$ for all $s\in H$.
\end{mydef}

\begin{myrem}
A Smarandache; RIPQ or LIPQ or IPQ(i.e a loop with at least a
non-trivial subquasigroup that is a RIPQ or LIPQ or IPQ) will now be
called a first Smarandache; RIPQ or LIPQ or
IPQ(S$_{1^{{\tiny\textrm{st}}}}$RIPQ or
S$_{1^{{\tiny\textrm{st}}}}$LIPQ or S$_{1^{{\tiny\textrm{st}}}}$IPQ
). It is easy to see that a S$_{2^{{\tiny\textrm{st}}}}$RIPQ or
S$_{2^{{\tiny\textrm{st}}}}$LIPQ or S$_{2^{{\tiny\textrm{st}}}}$IPQ
is a S$_{1^{{\tiny\textrm{st}}}}$RIPQ or
S$_{1^{{\tiny\textrm{st}}}}$LIPQ or S$_{1^{{\tiny\textrm{st}}}}$IPQ
respectively. But the reverse is not generally true.
\end{myrem}

\begin{mydef}\label{1:3}
Let $(G,\cdot )$ be a quasigroup(loop). The set $SYM(G, \cdot
)=SYM(G)$ of all bijections in $G$ forms a group called the
permutation(symmetric) group of $G$. The triple $(U, V, W)$ such
that $U, V, W\in SYM(G, \cdot )$ is called an autotopism of $G$ if
and only if
\begin{displaymath}
xU\cdot yV=(x\cdot y)W~\forall ~x, y\in G.
\end{displaymath}
The group of autotopisms of $G$ is denoted by $AUT(G, \cdot
)=AUT(G)$.

Let $(G_H,\cdot )$ be a special quasigroup(loop). The set $SSYM(G_H,
\cdot )=SSYM(G_H)$ of all Smarandache bijections(S-bijections) in
$G_H$ i.e $A\in SYM(G_H)$ such that $A~:~H\to H$ forms a group
called the Smarandache permutation(symmetric) group[S-permutation
group] of $G_H$. The triple $(U, V, W)$ such that $U, V, W\in
SSYM(G_H, \cdot )$ is called a first Smarandache
autotopism(S$_{1^{{\tiny\textrm{st}}}}$ autotopism) of $G_H$ if and
only if
\begin{displaymath}
xU\cdot yV=(x\cdot y)W~\forall ~x, y\in G_H.
\end{displaymath}
If their set forms a group under componentwise multiplication, it is
called the first Smarandache autotopism
group(S$_{1^{{\tiny\textrm{st}}}}$ autotopism group) of $G_H$ and is
denoted by S$_{1^{{\tiny\textrm{st}}}}AUT(G_H,\cdot
)=\textrm{S}_{1^{{\tiny\textrm{st}}}}AUT(G_H)$.

The triple $(U, V, W)$ such that $U,W\in SYM(G, \cdot )$ and $V\in
SSYM(G_H, \cdot )$ is called a second right Smarandache
autotopism(S$_{2^{{\tiny\textrm{nd}}}}$ right autotopism) of $G_H$
if and only if
\begin{displaymath}
xU\cdot sV=(x\cdot s)W~\forall~x\in G~\textrm{and}~s\in H.
\end{displaymath}
If their set forms a group under componentwise multiplication, it is
called the second right Smarandache autotopism
group(S$_{2^{{\tiny\textrm{nd}}}}$ right autotopism group) of $G_H$
and is denoted by S$_{2^{{\tiny\textrm{nd}}}}RAUT(G_H,\cdot
)=\textrm{S}_{2^{{\tiny\textrm{nd}}}}RAUT(G_H)$.

The triple $(U, V, W)$ such that $V,W\in SYM(G, \cdot )$ and $U\in
SSYM(G_H, \cdot )$  is called a second left Smarandache
autotopism(S$_{2^{{\tiny\textrm{nd}}}}$ left autotopism) of $G_H$ if
and only if
\begin{displaymath}
sU\cdot yV=(s\cdot y)W~\forall~y\in G~\textrm{and}~s\in
H.\end{displaymath} If their set forms a group under componentwise
multiplication, it is called the second left Smarandache autotopism
group(S$_{2^{{\tiny\textrm{nd}}}}$ left autotopism group) of $G_H$
and is denoted by S$_{2^{{\tiny\textrm{nd}}}}LAUT(G_H,\cdot
)=\textrm{S}_{2^{{\tiny\textrm{nd}}}}LAUT(G_H)$.

Let $(G_H,\cdot )$ be a special quasigroup(loop) with identity
element $e$. A mapping $T\in SSYM(G_H)$ is called a first
Smarandache semi-automorphism(S$_{1^{{\tiny\textrm{st}}}}$
semi-automorphism) if and only if $eT=e$ and
\begin{displaymath}
(xy\cdot x)T=(xT\cdot yT)xT~\textrm{for all}~x,y\in G.
\end{displaymath}

A mapping $T\in SSYM(G_H)$ is called a second Smarandache
semi-automorphism(S$_{2^{{\tiny\textrm{nd}}}}$ semi-automorphism) if
and only if $eT=e$ and
\begin{displaymath}
(sy\cdot s)T=(sT\cdot yT)sT~\textrm{for all}~y\in G~\textrm{and
all}~s\in H.
\end{displaymath}

A special loop $(G_H,\cdot )$ is called a first Smarandache
semi-automorphism inverse property
loop(S$_{1^{{\tiny\textrm{st}}}}$SAIPL) if and only if $J_\rho$ is a
S$_{1^{{\tiny\textrm{st}}}}$ semi-automorphism.

A special loop $(G_H,\cdot )$ is called a second Smarandache
semi-automorphism inverse property
loop(S$_{2^{{\tiny\textrm{nd}}}}$SAIPL) if and only if $J_\rho$ is a
S$_{2^{{\tiny\textrm{nd}}}}$ semi-automorphism.

Let $(G_H,\cdot )$ be a special quasigroup(loop). A mapping $A\in
SSYM(G_H)$ is a
\begin{enumerate}
\item first Smarandache pseudo-automorphism(S$_{1^{{\tiny\textrm{st}}}}$ pseudo-automorphism) of $G_H$ if and only if
there exists a $c\in H$ such that
$(A,AR_c,AR_c)\in\textrm{S}_{1^{{\tiny\textrm{st}}}}AUT(G_H)$. $c$
is reffered to as the first Smarandache
companion(S$_{1^{{\tiny\textrm{st}}}}$ companion) of $A$. The set of
such $A$s' is denoted by S$_{1^{{\tiny\textrm{st}}}}PAUT(G_H,\cdot
)=\textrm{S}_{1^{{\tiny\textrm{st}}}}PAUT(G_H)$.
\item second right Smarandache pseudo-automorphism(S$_{2^{{\tiny\textrm{nd}}}}$ right pseudo-automorphism) of $G_H$ if and only if
there exists a $c\in H$ such that
$(A,AR_c,AR_c)\in\textrm{S}_{2^{{\tiny\textrm{nd}}}}RAUT(G_H)$. $c$
is reffered to as the second right Smarandache
companion(S$_{2^{{\tiny\textrm{nd}}}}$ right companion) of $A$. The
set of such $A$s' is denoted by
S$_{2^{{\tiny\textrm{nd}}}}RPAUT(G_H,\cdot
)=\textrm{S}_{2^{{\tiny\textrm{nd}}}}RPAUT(G_H)$.
\item second left Smarandache pseudo-automorphism(S$_{2^{{\tiny\textrm{nd}}}}$ left pseudo-automorphism) of $G_H$ if and only if
there exists a $c\in H$ such that
$(A,AR_c,AR_c)\in\textrm{S}_{2^{{\tiny\textrm{nd}}}}LAUT(G_H)$. $c$
is reffered to as the second left Smarandache
companion(S$_{2^{{\tiny\textrm{nd}}}}$ left companion) of $A$. The
set of such $A$s' is denoted by
S$_{2^{{\tiny\textrm{nd}}}}LPAUT(G_H,\cdot
)=\textrm{S}_{2^{{\tiny\textrm{nd}}}}LPAUT(G_H)$.
\end{enumerate}
\end{mydef}

\section{Main Results}
\begin{myth}\label{1:4}
Let the special loop $(G_H,\cdot )$ be a
S$_{2^{{\tiny\textrm{nd}}}}$BL. Then it is both a
S$_{2^{{\tiny\textrm{nd}}}}$RIPL and a
S$_{2^{{\tiny\textrm{nd}}}}$RAPL.
\end{myth}
{\bf Proof}\\
\begin{enumerate}
\item In the S$_{2^{{\tiny\textrm{nd}}}}$BI, substitute $z=s^\rho$, then
$(xs\cdot s^\rho)s=x(ss^\rho\cdot s)=xs$ for all $x\in G$ and $s\in
H$. Hence, $xs\cdot s^\rho=x$ which is the
S$_{2^{{\tiny\textrm{nd}}}}$RIP.
\item In the S$_{2^{{\tiny\textrm{nd}}}}$BI, substitute $z=e$ and get
$xs\cdot s=x\cdot ss$ for all $x\in G$ and $s\in H$. Which is the
S$_{2^{{\tiny\textrm{nd}}}}$RAP.
\end{enumerate}

\begin{myrem}
Following Theorem~\ref{1:4}, we know that if a special loop
$(G_H,\cdot )$ is a S$_{2^{{\tiny\textrm{nd}}}}$BL, then its special
subloop $(H,\cdot )$ is a Bol loop. Hence, $s^{-1}=s^\lambda
=s^\rho$ for all $s\in H$. So, if $n\in\mathbb{Z}^+$, define $xs^n$
recursively by $s^0=e$ and $s^n=s^{n-1}\cdot s$. For any $n\in
\mathbb{Z}^-$, define $s^n$ by $s^n=(s^{-1})^{|n|}$.
\end{myrem}

\begin{myth}\label{1:5}
If $(G_H,\cdot )$ is a S$_{2^{{\tiny\textrm{nd}}}}$BL, then
\begin{equation}\label{eq:2}
xs^n=xs^{n-1}\cdot s=xs\cdot s^{n-1}
\end{equation}
for all $x\in G$, $s\in H$ and $n\in \mathbb{Z}$.
\end{myth}
{\bf Proof}\\
Trivialy, (\ref{eq:2}) holds for $n=0$ and $n=1$. Now assume for
$k>1$,
\begin{equation}\label{eq:3}
xs^k=xs^{k-1}\cdot s=xs\cdot s^{k-1}
\end{equation}
for all $x\in G$, $s\in H$. In particular, $s^k=s^{k-1}\cdot
s=s\cdot s^{k-1}$ for all $s\in H$. So, $xs^{k+1}=x\cdot
s^ks=x(ss^{k-1}\cdot s)=(xs\cdot s^{k-1})s=xs^k\cdot s$ for all
$x\in G$, $s\in H$. Then, replacing $x$ by $xs$ in (\ref{eq:3}),
$xs\cdot s^k=(xs\cdot s^{k-1})s=x(ss^{k-1}\cdot s)=x(s^{k-1}s\cdot
s)=x\cdot s^ks=xs^{k+1}$ for all $x\in G$, $s\in H$.(Note that the
S$_{2^{{\tiny\textrm{nd}}}}$BI has been used twice.)

Thus, (\ref{eq:2}) holds for all integers $n\ge 0$.

Now, for all integers $n>0$ and all $x\in G$, $s\in H$, applying
(\ref{eq:2}) to $x$ and $s^{-1}$ gives
$x(s^{-1})^{n+1}=x(s^{-1})^n\cdot s^{-1}=xs^{-n}\cdot s^{-1}$, and
(\ref{eq:2}) applied to $xs$ and $s^{-1}$ gives $xs\cdot
(s^{-1})^{n+1}=(xs\cdot s^{-1})(s^{-1})^n=xs^{-n}$. Hence,
$xs^{-n}=xs^{-n-1}\cdot s=xs\cdot s^{-n-1}$ and the proof is
complete.(Note that the S$_{2^{{\tiny\textrm{nd}}}}$RIP of
Theorem~\ref{1:4} has been used.)

\begin{myth}\label{1:6}
If $(G_H,\cdot )$ is a S$_{2^{{\tiny\textrm{nd}}}}$BL, then
\begin{equation}\label{eq:4}
xs^m\cdot s^n=xs^{m+n}
\end{equation}
for all $x\in G$, $s\in H$ and $m,n\in \mathbb{Z}$.
\end{myth}
{\bf Proof}\\
The desired result clearly holds for $n=0$ and by Theorem~\ref{1:5},
it also holds for $n=1$.

For any integer $n>1$, assume that (\ref{eq:4}) holds for all $m\in
\mathbb{Z}$ and all $x\in G$, $s\in H$. Then, using
Theorem~\ref{1:5}, $xs^{m+n+1}=xs^{m+n}\cdot s=(xs^m\cdot
s^n)s=xs^m\cdot s^{n+1}$ for all $x\in G$, $s\in H$ and $m\in
\mathbb{Z}$. So, (\ref{eq:4}) holds for all $m\in \mathbb{Z}$ and
$n\in \mathbb{Z}^+$. Recall that $(s^n)^{-1}=s^{-n}$ for all $n\in
\mathbb{Z}^+$ and $s\in H$. Replacing $m$ by $m-n$, $xs^{m-n}\cdot
s^n=xs^m$ and, hence, $xs^{m-n}=xs^m\cdot (s^n)^{-1}=xs^m\cdot
s^{-n}$ for all $m\in \mathbb{Z}$ and $x\in G$, $s\in H$.

\begin{mycor}\label{1:7}
Every S$_{2^{{\tiny\textrm{nd}}}}$BL is a SRPAPL.
\end{mycor}
{\bf Proof}\\
When $n=1$, the SRPAP is true. When $n=2$, the SRPAP is the SRAP.
Let the SRPAP be true for $k\in \mathbb{Z}^+$; $R_{s^k}=R_s^k$ for
all $s\in H$. Then, by Theorem~\ref{1:6},
$R_s^{k+1}=R_s^kR_s=R_{s^k}R_s=R_{s^{k+1}}$ for all $s\in H$.

\begin{mylem}\label{1:8}
Let $(G_H,\cdot )$ be a special loop. Then,
S$_{1^{{\tiny\textrm{st}}}}AUT(G_H,\cdot )\le AUT(G_H,\cdot
)$,~S$_{2^{{\tiny\textrm{nd}}}}RAUT(G_H,\cdot )\le AUT(H,\cdot )$
and S$_{2^{{\tiny\textrm{nd}}}}LAUT(G_H,\cdot )\le AUT(H,\cdot )$.
But, S$_{2^{{\tiny\textrm{nd}}}}RAUT(G_H,\cdot )\not\le
AUT(G_H,\cdot )$ and S$_{2^{{\tiny\textrm{nd}}}}LAUT(G_H,\cdot
)\not\le AUT(G_H,\cdot )$.
\end{mylem}
{\bf Proof}\\
These are easily proved by using the definitions of the sets
relative to componentwise multiplication.

\begin{mylem}\label{1:9}
Let $(G_H,\cdot )$ be a special loop. Then,
S$_{2^{{\tiny\textrm{nd}}}}RAUT(G_H,\cdot )$ and
S$_{2^{{\tiny\textrm{nd}}}}LAUT(G_H,\cdot )$ are groups under
componentwise multiplication.
\end{mylem}
{\bf Proof}\\
These are easily proved by using the definitions of the sets
relative to componentwise multiplication.

\begin{mylem}\label{1:10}
Let $(G_H,\cdot )$ be a special loop.
\begin{enumerate}
\item If $(U, V,W)\in\textrm{S}_{2^{{\tiny\textrm{nd}}}}RAUT(G_H,\cdot )$ and
$G_H$ has the S$_{2^{{\tiny\textrm{nd}}}}$RIP, then $(W, J_\rho
VJ_\rho ,U)\in\textrm{S}_{2^{{\tiny\textrm{nd}}}}RAUT(G_H,\cdot )$.
\item If $(U, V,W)\in\textrm{S}_{2^{{\tiny\textrm{nd}}}}LAUT(G_H,\cdot )$ and
$G_H$ has the S$_{2^{{\tiny\textrm{nd}}}}$LIP, then $(J_\lambda U,
W,V)\in\textrm{S}_{2^{{\tiny\textrm{nd}}}}LAUT(G_H,\cdot )$.
\end{enumerate}
\end{mylem}
{\bf Proof}\\
\begin{enumerate}
\item $(U, V,W)\in\textrm{S}_{2^{{\tiny\textrm{nd}}}}RAUT(G_H,\cdot )$ implies
that $xU\cdot sV=(x\cdot s)W$ for all $x\in G$ and $s\in H$. So,
$(xU\cdot sV)(sV)^\rho=(x\cdot s)W\cdot (sV)^\rho\Rightarrow
xU=(xs^\rho )W\cdot (s^\rho V)^\rho\Rightarrow (xs)U=(xs\cdot s^\rho
)W\cdot (s^\rho V)^\rho\Rightarrow (xs)U=xW\cdot sJ_\rho
VJ_\rho\Rightarrow (W, J_\rho VJ_\rho
,U)\in\textrm{S}_{2^{{\tiny\textrm{nd}}}}RAUT(G_H,\cdot )$.
\item $(U, V,W)\in\textrm{S}_{2^{{\tiny\textrm{nd}}}}LAUT(G_H,\cdot )$ implies
that $sU\cdot xV=(s\cdot x)W$ for all $x\in G$ and $s\in H$. So,
$(sU)^\lambda\cdot (sU\cdot xV)=(sU)^\lambda\cdot (s\cdot
x)W\Rightarrow xV=(sU)^\lambda\cdot (sx)W\Rightarrow xV=(s^\lambda
U)^\lambda\cdot (s^\lambda x)W\Rightarrow (sx)V=(s^\lambda
U)^\lambda\cdot (s^\lambda\cdot sx)W\Rightarrow (sx)V=sJ_\lambda
UJ_\lambda\cdot xW\Rightarrow (J_\lambda U,
W,V)\in\textrm{S}_{2^{{\tiny\textrm{nd}}}}LAUT(G_H,\cdot )$.
\end{enumerate}

\begin{myth}\label{1:11}
Let $(G_H,\cdot )$ be a special loop. $(G_H,\cdot )$ is a
S$_{2^{{\tiny\textrm{nd}}}}$BL if and only if
$(R_s^{-1},L_sR_s,R_s)\in\textrm{S}_{1^{{\tiny\textrm{st}}}}AUT(G_H,\cdot
)$.
\end{myth}
{\bf Proof}\\
$G_H$ is a S$_{2^{{\tiny\textrm{nd}}}}$BL iff $(xs\cdot
z)s=x(sz\cdot s)$ for all $x,z\in G$ and $s\in H$ iff $(xR_s\cdot
z)R_s=x(zL_sR_s)$ iff $(xz)R_s=xR_s^{-1}\cdot zL_sR_s$ iff
$(R_s^{-1},L_sR_s,R_s)\in\textrm{S}_{1^{{\tiny\textrm{st}}}}AUT(G_H,\cdot
)$.

\begin{myth}\label{1:12}
Let $(G_H,\cdot )$ be a S$_{2^{{\tiny\textrm{nd}}}}$BL. $G_H$ is a
S$_{2^{{\tiny\textrm{nd}}}}$SAIPL if and only if $G_H$ is a
S$_{3^{{\tiny\textrm{rd}}}}$RIPL.
\end{myth}
{\bf Proof}\\
Keeping the S$_{2^{{\tiny\textrm{nd}}}}$BI and the
S$_{2^{{\tiny\textrm{nd}}}}$RIP in mind, it will be observed that if
$G_H$ is a S$_{3^{{\tiny\textrm{rd}}}}$RIPL, then $(sy\cdot
s)(s^\rho y^\rho\cdot s^\rho )=[((sy\cdot s)s^\rho )y^\rho]s^\rho
=(sy\cdot y^\rho)s^\rho =ss^\rho =e$. So, $(sy\cdot s)^\rho =s^\rho
y^\rho\cdot s^\rho$. The proof of the necessary part follows by the
reverse process.

\begin{myth}\label{1:13}
Let $(G_H,\cdot )$ be a S$_{2^{{\tiny\textrm{nd}}}}$BL. If $(U,
T,U)\in\textrm{S}_{1^{{\tiny\textrm{st}}}}AUT(G_H,\cdot )$. Then,
$T$ is a S$_{2^{{\tiny\textrm{nd}}}}$ semi-automorphism.
\end{myth}
{\bf Proof}\\
If $(U, T,U)\in\textrm{S}_{1^{{\tiny\textrm{st}}}}AUT(G_H,\cdot )$,
then, $(U, T,U)\in\textrm{S}_{2^{{\tiny\textrm{nd}}}}RAUT(G_H,\cdot
)\cap\textrm{S}_{2^{{\tiny\textrm{nd}}}}LAUT(G_H,\cdot )$.

Let $(U, T,U)\in\textrm{S}_{2^{{\tiny\textrm{nd}}}}RAUT(G_H,\cdot
)$, then $xU\cdot sT=(xs)U$ for all $x\in G$ and $s\in H$. Set
$s=e$, then $eT=e$. Let $u=eU$, then $u\in H$ since $(U,
T,U)\in\textrm{S}_{2^{{\tiny\textrm{nd}}}}LAUT(G_H,\cdot )$. For
$x=e$, $U=TL_u$. So, $xTL_u\cdot sT=(xs)TL_u$ for all $x\in G$ and
$s\in H$. Thus,
\begin{equation}\label{eq:5}
(u\cdot xT)\cdot sT=u\cdot (xs)T.
\end{equation}
Replace $x$ by $sx$ in (\ref{eq:5}), to get
\begin{equation}\label{eq:6}
[u\cdot (sx)T]\cdot sT=u\cdot (sx\cdot s)T.
\end{equation}
$(U, T,U)\in\textrm{S}_{2^{{\tiny\textrm{nd}}}}LAUT(G_H,\cdot )$
implies that $sU\cdot xT=(sx)U$ for all $x\in G$ and $s\in H$
implies $sTL_u\cdot xT=(sx)TL_u$ implies $(u\cdot sT)\cdot xT=u\cdot
(sx)T$. Using this in (\ref{eq:6}) gives $[(u\cdot sT)\cdot xT]\cdot
sT=u\cdot (sx\cdot s)T$. By the S$_{2^{{\tiny\textrm{nd}}}}$BI,
$u[(sT\cdot xT)\cdot sT]=u\cdot (sx\cdot s)T\Rightarrow (sT\cdot
xT)\cdot sT=(sx\cdot s)T$.

\begin{mycor}\label{1:14}
Let $(G_H,\cdot )$ be a S$_{2^{{\tiny\textrm{nd}}}}$BL that is a
Smarandache right nuclear square special loop. Then, $L_sR_s^{-1}$
is a S$_{2^{{\tiny\textrm{nd}}}}$ semi-automorphism.
\end{mycor}
{\bf Proof}\\
$s^2\in SN_\rho (G_H)$ for all $s\in H$ iff $xy\cdot s^2=x\cdot
ys^2$ iff $(xy)R_{s^2}=x\cdot yR_{s^2}$ iff $(xy)R_s^2=x\cdot
yR_s^2$($\because$ of S$_{2^{{\tiny\textrm{nd}}}}$RAP) iff
$(I,R_s^2,R_s^2)\in \textrm{S}_{1^{{\tiny\textrm{st}}}}AUT(G_H,\cdot
)$ iff $(I,R_s^{-2},R_s^{-2})\in
\textrm{S}_{1^{{\tiny\textrm{st}}}}AUT(G_H,\cdot )$. Recall from
Theorem~\ref{1:11} that,
$(R_s^{-1},L_sR_s,R_s)\in\textrm{S}_{1^{{\tiny\textrm{st}}}}AUT(G_H,\cdot
)$. So,
$(R_s^{-1},L_sR_s,R_s)(I,R_s^{-2},R_s^{-2})=(R_s^{-1},L_sR_s^{-1},R_s^{-1})\in\textrm{S}_{1^{{\tiny\textrm{st}}}}AUT(G_H,\cdot
)\Rightarrow L_sR_s^{-1}$ is a S$_{2^{{\tiny\textrm{nd}}}}$
semi-automorphism by Theorem~\ref{1:13}.

\begin{mycor}\label{1:15}
If a S$_{2^{{\tiny\textrm{nd}}}}$BL is of Smarandache exponent $2$,
then, $L_sR_s^{-1}$ is a S$_{2^{{\tiny\textrm{nd}}}}$
semi-automorphism.
\end{mycor}
{\bf Proof}\\
These follows from Theorem~\ref{1:14}.

\begin{myth}\label{1:16}
Let $(G_H,\cdot )$ be a S$_{2^{{\tiny\textrm{nd}}}}$BL. Let $(U,
V,W)\in\textrm{S}_{1^{{\tiny\textrm{st}}}}AUT(G_H,\cdot )$, $s_1=eU$
and $s_2=eV$. Then, $A=UR_s^{-1}\in
\textrm{S}_{1^{{\tiny\textrm{st}}}}PAUT(G_H)$ with
S$_{1^{{\tiny\textrm{st}}}}$ companion $c=s_1s_2\cdot s_1$ such that
$(U, V,W)=(A,AR_c,AR_c)(R_s^{-1},L_sR_s,R_s)^{-1}$.
\end{myth}
{\bf Proof}\\
By Theorem~\ref{1:11},
$(R_s^{-1},L_sR_s,R_s)\in\textrm{S}_{1^{{\tiny\textrm{st}}}}AUT(G_H,\cdot
)$ for all $s\in H$. Hence, $(A,B,C)=(U,
V,W)(R_{s_1}^{-1},L_{s_1}R_{s_1},R_{s_1})=(UR_{s_1}^{-1},VL_{s_1}R_{s_1},WR_{s_1})\in\textrm{S}_{1^{{\tiny\textrm{st}}}}AUT(G_H,\cdot
)\Rightarrow~A=UR_{s_1}^{-1},~B=VL_{s_1}R_{s_1}$ and $C=WR_{s_1}$.
That is, $aA\cdot bB=(ab)C$ for all $a,b\in G_H$. Since $eA=e$, then
setting $a=e$, $B=C$. Then for $b=e$, $B=AR_{eB}$. But
$eB=eVL_{s_1}R_{s_1}=s_1s_2\cdot s_1$. Thus,
$(A,AR_{eB},AR_{eB})\in\textrm{S}_{1^{{\tiny\textrm{st}}}}AUT(G_H,\cdot
)\Rightarrow~A\in\textrm{S}_{1^{{\tiny\textrm{st}}}}PAUT(G_H,\cdot
)$ with S$_{1^{{\tiny\textrm{st}}}}$ companion $c=s_1s_2\cdot s_1\in
H$.

\begin{myth}\label{1:17}
Let $(G_H,\cdot )$ be a S$_{2^{{\tiny\textrm{nd}}}}$BL. Let $(U,
V,W)\in\textrm{S}_{2^{{\tiny\textrm{nd}}}}LAUT(G_H,\cdot
)\cap\textrm{S}_{2^{{\tiny\textrm{nd}}}}RAUT(G_H,\cdot )$, $s_1=eU$
and $s_2=eV$. Then, $A=UR_s^{-1}\in
\textrm{S}_{2^{{\tiny\textrm{nd}}}}LPAUT(G_H)\cap\textrm{S}_{2^{{\tiny\textrm{nd}}}}RPAUT(G_H)$
with S$_{2^{{\tiny\textrm{nd}}}}$ left companion and
S$_{2^{{\tiny\textrm{nd}}}}$ right companion $c=s_1s_2\cdot s_1$
such that $(U, V,W)=(A,AR_c,AR_c)(R_s^{-1},L_sR_s,R_s)^{-1}$.
\end{myth}
{\bf Proof}\\
The proof of this is very similar to the proof of
Theorem~\ref{1:16}.

\begin{myrem}
Every Bol loop is a S$_{2^{{\tiny\textrm{nd}}}}$BL. Most of the
results on basic properties of Bol loops in chapter~2 of
\cite{phd85} can easily be deduced from the results in this paper by
simply forcing $H$ to be equal to $G$.
\end{myrem}

\begin{myquest}
Let $(G_H,\cdot )$ be a special quasigroup(loop). Are the sets
S$_{1^{{\tiny\textrm{st}}}}PAUT(G_H)$,
S$_{2^{{\tiny\textrm{nd}}}}RPAUT(G_H)$ and
S$_{2^{{\tiny\textrm{nd}}}}LPAUT(G_H)$ groups under mapping
composition?
\end{myquest}
\begin{myquest}
Let $(G_H,\cdot )$ be a special quasigroup(loop). Can we find a
general method(i.e not an "acceptable"
S$_{2^{{\tiny\textrm{nd}}}}$BL with carrier set $\mathbb{N}$) of
constructing a S$_{2^{{\tiny\textrm{nd}}}}$BL that is not a Bol loop
just like Robinson \cite{phd85}, Solarin and Sharma \cite{phd5} were
able to use general methods to construct Bol loops.
\end{myquest}

\end{document}